\documentclass[a4paper,12pt]{amsart}
\usepackage{amsmath,inputenc,euscript,amssymb}
\usepackage{hyperref}
\title[NLS on hyperbolic space]{The nonlinear Schr\"odinger equation\\ on hyperbolic space}
\author{V. Banica}

\newtheorem{lemma}{Lemma}[section]
\newtheorem{theorem}[lemma]{Theorem}
\newtheorem{prop}[lemma]{Proposition}
\newtheorem{corollary}[lemma]{Corollary}

\newtheorem{remark}[lemma]{Remark}

\newcommand{\mo}[1]{|#1|}

\newcommand{\RR}{\mathbb{R}}
\newcommand{\TT}{\mathbb{T}}
\newcommand{\HH}{\mathbb{H}}

\newcommand{\DP}[1]{{\partial_ #1}}

\def\tend{{\rightarrow}}
\newcommand{\R}{\mathbb{R}}

\thanks{Research supported by HARP (HPRN-CT-2001-00273) and 
HYKE (HPRN-CT-2002-00282) research projects. Email address : Valeria.Banica at univ-evry.fr}
\begin{document}

\begin{abstract}
In this article we study some aspects of dispersive and concentration 
phenomena for the Schr\"odinger equation posed on hyperbolic space 
$\mathbb{H}^n$, in order to see if the negative curvature of the manifold 
gets the dynamics more stable than in the Euclidean case. 
It is indeed the case for the dispersive properties : we prove that the dispersion inequality is 
valid, in a stronger form than the one on $\mathbb{R}^n$. 
However, the geometry does not have enough of an effect to avoid the 
concentration phenomena and the picture is actually worse than expected. The critical nonlinearity power for blow-up turns out to be the same as in the euclidean case, and we prove that there are more explosive solutions for critical and supercritical nonlinearities. 

\end{abstract}
\maketitle
\section{Introduction}
Let $(M,g)$ be a Riemannian manifold and let $\Delta_M$ be the 
Laplace-Beltrami operator. 
The study of the nonlinear Schr{\"o}dinger equation
$$\left\{\begin{array}{c}
i\DP{t}u+\Delta_M u=V'(\mo{u}^2)u,\\
u(0)=u_0,
\end{array}\right.$$
where $u$ is a space-time function with complex values, and $V$ is a real 
function with controlled growth at infinity, was motivated by number of 
problems coming from Physics. 

It is known that the geometry influences the dynamics of the equation. 
Instability phenomena appear, even in the defocusing case.\\
For instance, on the one hand, for the cubic defocusing Schr\"odinger equation on the sphere
$\mathbb{S}^2$
$$\left\{\begin{array}{c}
i\partial_t u+\Delta_{\mathbb{S}^2}u=\vert u\vert^2u,\\
u(0,x)\in H^s(\mathbb{S}^2)
\end{array}\right.,$$
the critical regularity index for having the uniform continuity of the flow 
on the bounded sets of $H^s$ is $s=\frac{1}{4}$, as it was proved recently by 
Burq-G\'erard-Tzvetkov (\cite{BGT12}, \cite{BGT32}, see also \cite{Val2}).\\
On the other hand Bourgain (\cite{Bo}) and Cazenave-Weissler (\cite{CW}) 
have proved that the Cauchy problem of the same equation, 
considered on $\TT^2$ and 
on $\RR^2$ respectively, is $H^\epsilon$ well-posed for all positive 
$\epsilon$, meanwhile for negative 
$s$ instability phenomena appear (\cite{BGT12},\cite{CCT}).
It follows that the critical regularity index for the flat torus and for 
$\RR^2$ is zero.

Hence these results point out the importance played by the geometry of the
manifold in the dynamics of the equation.

It is expected that the positive curvature generates the 
differences, since for having the instability result on the sphere the 
dynamics of spherical harmonic concentrated on closed geodesics are studied.

In this article we study dispersive and concentration phenomena for the 
Schr\"odinger equation posed on hyperbolic space $\mathbb{H}^n$,
manifold of negative curvature, 
expecting that the dynamics are stabler than in the Euclidean case.\\

We shall define in \S 2 the hyperbolic space and the tools used on
it. For the moment, let us introduce some notations. We denote by $0$ the
origin of the hyperbolic space, $0=(1,0,..,0)$. In the sequel we shall
use the $L^p$ spaces on $\mathbb{H}^n$
$$L^p=L^p(\mathbb{H}^n)=L^p(d\Omega),$$
and weighted spaces, defined by
$$L^p(w)=L^p(w\,d\Omega),$$
where $d\Omega$ denotes the measure on hyperbolic space
$\mathbb{H}^n$ and $w$ is a function on $\mathbb{H}^n$.\\

First, we shall treat the linear equation. We obtain an explicit
representation of the solutions.

\begin{theorem}\label{exactsol}
The solution of the linear Schr\"odinger equation posed on 
hyperbolic space $\mathbb{H}^n$ 
$$(SL)\left\{\begin{array}{c}
i\partial_t u+\Delta_{\mathbb{H}^n}u=0\\
u(0,x)=u_0
\end{array}\right.,$$
is, for $n\geq 3$ odd,
\begin{equation}\label{odd}
u(t,\Omega)=\frac{c}{|t|^\frac{1}{2}}e^{-it\frac{(n-1)^2}{4}}\int_{\mathbb{H}^n}u_o(\Omega')\left(\frac{\partial_\rho}{\sinh\rho}\right)^\frac{n-1}{2}e^{i\frac{\rho^2}{4t}}\,d\Omega',\end{equation}
and for $n\geq 2$ even
\begin{equation}\label{even}
u(t,\Omega)=\frac{c}{|t|^\frac{3}{2}}e^{-it\frac{(n-1)^2}{4}}\int_{\mathbb{H}^n}u_o(\Omega')\left(\frac{\partial_\rho}{\sinh
\rho}\right)^\frac{n-2}{2}\int_\rho^\infty\frac{e^{i\frac{s^2}{4t}}s}{\sqrt{\cosh
s-\cosh\rho}}\,ds\,d\Omega',\end{equation}
where by $\rho$ we denote the hyperbolic distance between the points
$\Omega$ and $\Omega'$.

\end{theorem}
As harmonic analysis can be done on $\mathbb{H}^n$, 
the proof is based on the representation of the solution via the Fourier 
transform and on calculus of oscillatory integrals. There are many
points in common with the proof of the inverse Fourier formula on hyperbolic space.

This explicit representation of the solution allows us to study the
dispersive properties and to obtain the following results.

\begin{theorem}\label{disp}

i) For all dimension $n\geq 2$, the solution satisfies the following local 
dispersion inequality 
\begin{equation}\label{disph}
|u(t,\Omega)|\leq
c\frac{1}{|t|^\frac{n}{2}}\int_{\mathbb{H}^n}|u_0(\Omega')|\left(\frac{\rho}{\sinh\rho}\right)^\frac{n-1}{2}\,d\Omega'.\end{equation}

For large times, the following dispersion estimate holds for all $n\geq 3$ odd,
\begin{equation}\label{disphtlarge}
|u(t,\Omega)|\leq
c\frac{1}{|t|^\frac{3}{2}}\int_{\mathbb{H}^n}|u_0(\Omega')|\left(\frac{\rho}{\sinh\rho}\right)^\frac{n-1}{2}\,d\Omega',\end{equation}
and for all $n\geq 2$ even
\begin{equation}\label{disphtlarge2}
|u(t,\Omega)|\leq
c\frac{1}{|t|^\frac{3}{2}}\int_{\mathbb{H}^n}|u_0(\Omega')|\left(\frac{\rho}{\sinh\rho}\right)^\frac{n-1}{2}\,\frac{1+\rho}{\sqrt{\rho}}\,d\Omega',\end{equation}

ii) 
Moreover, for $n\geq 3$ and radial initial data, we have a weighted-space local dispersion inequality
\begin{equation}\label{dispw}
\|u(t)\|_{L^\infty(w)}\leq \frac{c}{|t|^\frac{n}{2}}
\|u_0\|_{L^1(w^{-1})},\end{equation}
where the weight is
$$w(\Omega)=\frac{\sinh d(0,\Omega)}{d(0,\Omega)}.$$
Finally, for $n\geq 3$, for a finite time $T$ and for radial initial data, we have the local Strichartz weighted estimates
\begin{equation}\label{Strw}
\|u\|_{L^p\left([0,T],\,
L^{q}\left(w^{q-2}\right)\right)}\leq 
c\|u_0\|_{L^2},\end{equation}
for all pairs $(p,q)$ satisfying $\frac{2}{p}+\frac{n}{q}=\frac{n}{2}$, $(p,q,n)\neq (2,\infty,2)$ and $2\leq p,q$.
\end{theorem}

The new term that appears in (\ref{disph}),
$$\frac{\rho}{\sinh \rho},$$
is specific to hyperbolic space ; it comes from the Harish-Chandra 
coefficient and from a kernel that involves the eigenfunctions of the Laplace-Beltrami operator. 
This new term informs us that, apart from the classical decay in time,
we have a new one, in space.

\begin{remark}
In dimension $3$, all the dispersion estimates from Theorem
\ref{disp} are global in time, which is not the case for other dimensions.
\end{remark}

\begin{remark}\label{dispsupp}
For dimensions allowed in Theorem \ref{disp}, away from the support of the initial data, the
decay is stronger than in the $\mathbb{R}^n$ case. More precisely, for initial data with support a domain $A$ of $\mathbb{H}^n$, the
solution satisfies, for any point $\Omega$ not included in $A$, at small times,
$$|u(t,\Omega)|\left(
\frac{\sinh d(\Omega,A)}{d(\Omega,A)}\right)^\frac{n-1}{2}\leq
\frac{c}{|t|^\frac{n}{2}}\int_{\mathbb{H}^n}|u_0(\Omega')|
\left(\frac{d(\Omega',^cA)}{\sinh d(\Omega',^cA)}\right)^\frac{n-1}{2}\,d\Omega'.$$
\end{remark}

\begin{remark}\label{disptime}
In dimension $3$, for radial initial data in a more restrictive space, that is
$L^1(\widetilde{w}^{-1})$ with the weight 
$$\widetilde{w}(\Omega)=\sinh d(\Omega,0),$$ we
obtain a local dispersion-type estimate stronger in time,
$$\|u(t)\|_{L^\infty(\widetilde{w})}\leq \frac{c}{|t|^{\frac{1}{2}}}
\|u_0\|_{L^1(\widetilde{w}^{-1})}.$$
\end{remark}

However, it is only in (\ref{dispw}), in the radial case, that we
obtained the dispersion stated in weighted spaces, independently of
the initial data.

Let us also notice that the improvements in the dispersive estimates
stated in Theorem \ref{disp} are better decays away from the
origin. Near the origin -where the Laplace-Beltrami operator on 
hyperbolic space is almost the one on $\mathbb{R}^n$- we have the same
decays as in the euclidean case. The infinite speed of propagation of
the Schr\"odinger operator is not sufficient to impose on the
dispersion estimate at the origin an influence from the metric at infinity.

As a consequence of (\ref{disph}), we obtain the classical local dispersion 
inequality
$$\|u(t)\|_{L^\infty}\leq \frac{c}{|t|^{\frac{n}{2}}}
\|u_0\|_{L^1}.$$\\
From the dispersion inequality, by using the TT* functional analysis 
argument (\cite{To}), the local Strichartz estimates are obtained
too. 
We recall here that for the wave equation on hyperbolic space, the Strichartz estimates have been proved recently by Tataru (\cite{Ta}).

Let us now remark that on the sphere, the local Strichartz-type 
estimates are known to hold with a loss of $\frac{1}{p}$ derivatives, 
that is the 
$L^pL^q$ norms of the solutions are controlled by the 
$H^\frac{1}{p}$ norm of the initial data, instead of the 
$L^2$ norm (\cite{BGT0}). This shows that on hyperbolic space, the
dispersion estimates are stronger.

\begin{remark}  An important problem related to the linear Schr\"odinger equation is the problem of the potential. The radial Schr\"odinger equation perturbed with a rough time dependent potential on hyperbolic space was recently treated by Pierfelice in (\cite{Vi}). By using Theorem \ref{disp} ii) and fixed point arguments she obtains for the radial perturbed Cauchy problem the Strichartz estimates with the same weight, pointing out again the influence of the negative curvature on the dispersive properties. These estimates hold even in the case of more general nonlinearities. 
\end{remark}

\begin{remark}
It is expected that the results obtained should hold on all symmetric spaces of
rank 1.

Moreover, let us notice that the dispersion kernel obtained is
close to the kernel of the heat operator $e^{t\Delta}$, which has 
been studied intensively on hyperbolic space, usually for $t$ real,
but also for complex $t$ with $\Re
t>0$ (\cite{D}, \cite{An}, \cite{DM}). In view of this work, we
expect that the dispersion estimate for large time can be improved by
an additional decay, for all $n\geq 2$ :
$$|u(t,\Omega)|\leq 
c\frac{1}{|t|^\frac{3}{2}}\int_{\mathbb{H}^n}|u_0(\Omega')|
\left(\frac{\rho}{\sinh\rho}\right)^\frac{n-1}{2}
\left(\frac{1}{1+\rho}\right)^\frac{n-3}{2}\,d\Omega'.$$
\end{remark}

\bigskip

Let us turn now our attention to the concentration phenomena for the
nonlinear equation, and more 
precisely to blow-up of solutions in the sense of the explosion in
finite time of the $L^2$ norm of their gradient.

We consider the focusing Schr\"odinger equation with power nonlinearity
$$(S)\left\{\begin{array}{c}
i\partial_t u+\Delta_{\mathbb{H}^n}u+\vert u\vert^{p-1}u=0\\
u(0,x)=u_0\in H^1(\mathbb{H}^n)
\end{array}\right..$$

First, let us notice that, as for the Euclidean case, since integration by parts works on hyperbolic space, the mass of a solution of $(S)$
$$\int_{\mathbb{H}^n}\mo{u}^2d\Omega,$$
and its energy 
$$E(u)=
\frac{1}{2}\int_{\mathbb{H}^n}\mo{\nabla u}^2d\Omega-
\frac{1}{p+1}\int_{\mathbb{H}^n}\mo{u}^{p+1}d\Omega,$$
are conserved in time.

By using the classical Strichartz estimates, derived from the
dispersion estimate (\ref{disph}) on the linear equation, 
one can obtain, as in the Euclidean case (\cite{GV},\cite{Ya}), the $H^1$ local existence of solutions of the equation $(S)$ with $p\leq 1+\frac{4}{n-2}$.

The argument used on $\mathbb{R}^n$ in the most of the blow-up results is the scale invariance. 
The construction of explicit blow-up solutions is done using the 
pseudo-conformal invariance. Informations on the blow-up solutions are 
obtained by studying the dilatations of the solution and by using the 
virial identity (\cite{Me}). 
The main difficulty when working on hyperbolic space is that these 
techniques and notions do not have an obvious analogue. 

One can expect that maybe there will be no blow-up solutions, due to the negative curvature of the manifold. 
In fact it turns out that the geometry does not have enough of an effect to avoid the 
concentration phenomena. Moreover, we have more types of initial data generating blow-up in finite time than in the Euclidean case. 
\begin{theorem}\label{blup}

For $p<1+\frac{4}{n}$, the solutions of the equation $(S)$ are global
in $H^1(\mathbb{H}^n)$. 
The global existence still holds for the power $p=1+\frac{4}{n}$ and
for initial data of mass smaller than a certain constant.\par
However, for $p\geq 1+\frac{4}{n}$, blow-up solutions exist. More
precisely, if the initial data is radial, of finite variance 
$$\int_{\mathbb{H}^n}|u_0(\Omega)|^2d^2(0,\Omega)d\Omega<+\infty,$$
and its energy satisfies
$$E(u_0)< c_n\|u_0\|^2_2,$$
then the solution blows up in finite time. 
Here $c_n$ is a geometric positive constant given by
$$c_n=\frac{\inf \Delta^2_{\mathbb{H}^n}d(0,\cdot)}{16}.$$

\end{theorem}
The proof of the global existence up to the critical power $p=1+\frac{4}{n}$ -which is also the one of the Euclidean space- involves 
the Sobolev 
embeddings on $\mathbb{H}^n$ (\cite{He}), together with the 
conservation laws of the equation. 
The existence of blow-up solutions is shown by
adapting the arguments of Glassey and Zakharov on $\mathbb{R}^n$
(\cite{Gl}, \cite{Za}) and of Kavian on a star-shaped domain of
$\mathbb{R}^n$ (\cite{Ka}), to the particular metric on $\mathbb{H}^n$.

\begin{remark}
We obtain a new virial identity adapted to hyperbolic space that 
allows us to
conclude that solutions of null energy and finite virial blow up in
finite time. 
Probably this new concentration phenomena is the consequence of the
improved dispersion of the
solutions. We recall that this is not the case on $\mathbb{R}^n$     
for the critical power,
since the ground state, that is the unique positive solution of the
elliptic equation
$$\Delta_{\mathbb{R}^n} Q+Q^p=Q,$$
gives us the global solution $e^{it}Q$ of null energy.\\ Moreover, on
$\mathbb{R}^n$, the mass of the ground state is the critical mass for
having explosive solutions in critical power.
In our case, the Theorem \ref{blup} implies that the ground state on hyperbolic space has positive
energy and mass smaller than the expected critical mass for blow-up on
$\mathbb{H}^n$.
\end{remark}

It should be interesting to see what is the minimal blow-up mass, and to get informations about the blow-up speed and about the profiles of the explosive solutions.\\

This article is organized as follows. In the second section we give the 
definition of hyperbolic spaces and of the tools used on it, namely the 
Laplace-Beltrami operator and the Fourier transform. 
In \S3.1 we recall the proof of the representation of the
solutions and the dispersion inequality in the 
Euclidean case, then we give in \S3.2 the proof of
Theorem \ref{exactsol}. 
In \S4 we prove the first point of Theorem \ref{disp}, concerning the
$L^1-L^\infty$ estimates. 
The section \S5 concerns the weighted dispersion inequalities of the second
point of Theorem \ref{disp}, and of the Remarks \ref{dispsupp} and 
\ref{disptime}. 
The last section contains the proof of Theorem \ref{blup}. In Appendix A is proved a technical proposition, crucial for
obtaining the dispersion estimates in \S4 .\\

Acknowledgments : I am very grateful to V. Georgiev, to P. G\'erard and to F. Ricci for having encouraged my research in this direction. I thank T. Duyckaerts and V. Pierfelice for having carefully read some parts of the article.


\section{hyperbolic space}
\subsection{Definition}

We shall use here the model of hyperbolic space given by the upper branch of the hyperboloid. 

We define the hyperbolic space as being the following surface of 
$\RR^{n+1}$, given by the parametrization :
$$\mathbb{H}^n=\{\Omega=(t,x)\in\RR^{n+1}, (t,x)=(\cosh r,\sinh r\,\omega), 
r\geq 0, \omega\in\mathbb{S}^{n-1}\}.$$

Let us introduce the inner product on $\RR^{n+1}$
\begin{equation}\label{prod}
[x,y]=x_0y_0-....-x_ny_n.\end{equation}

An alternative definition for hyperbolic space is 
$$\mathbb{H}^n=\{x=(x_0,x_1,...,x_n)\in\mathbb{R}^{n+1} , [x,x]=1 , x_0>0\}.$$
This space is invariant under $SO(1,n)$, the group of Lorentz 
transformations of $\RR^{n+1}$ that preserve this inner product. 

One has 
$$dt=\sinh r dr, dx=\cosh r\,\omega \,dr+\sinh r \,d\omega,$$
and the metric induced on $\mathbb{H}^n$ by the Lorenzian metric on $\RR^{n+1}$
$$dl^2=-dt^2+dx^2,$$
is 
$$ds^2=dr^2+\sinh^2 r\,d\omega^2,$$
where $d\omega^2$ is the metric on the sphere $\mathbb{S}^{n-1}$.

The volume element is
$$\int_{\mathbb{H}^n} f(\Omega)\,d\Omega=\int_0^\infty\int_{\mathbb{S}^{n-1}}f(r,\omega)\sinh^{n-1}r\,dr\,d\omega.$$

The length of a curve 
$$\gamma (t)=(\cosh r(t),\sinh r(t) \,\omega(t)),$$ with $t$ varying from $a$ to 
$b$, is defined to be as usual,
$$L(\gamma)=\int_a^b\sqrt{|r'(t)|^2+|\sinh r(t)|^2|\omega'(t)|^2}dt.$$
The distance between two points of $X$ will be the infimum of the lengths of 
the paths connecting the points. 
By direct calculus, one has that the distance of a point to the origin of the hyperboloid $O=(1,0,...,0)$ is 
$$d((\cosh r,\sinh r\,\omega), O)=r.$$
More generally, the distance between two arbitrary points is
\begin{equation}\label{dist}
d(\Omega,\Omega')=\cosh ^{-1}([\Omega,\Omega']).\end{equation}

Starting from the general definition of the Laplace-Beltrami 
operator, one can find its expression on $\mathbb{H}^n$
$$\Delta_{\mathbb{H}^n}=\partial_r^2+(n-1)\frac{\cosh r}{\sinh r}\,\partial_r+
\frac{1}{\sinh^2r}\,\Delta_{\mathbb{S}^{n-1}}.$$
In Section $6$ we will denote $\Delta=\Delta_{\mathbb{H}^n}$.

\subsection{The Fourier transform}\par
$$ $$
\par
For $\theta\in\mathbb{S}^{n-1}$ and $\lambda$ a real number, the functions of the 
type
\begin{equation}\label{fp}
h_{\lambda,\theta}(\Omega)=[\Omega,\Lambda(\theta)]^{i\lambda-\frac{n-1}{2}},
\end{equation}
are generalized eigenfunctions of the Laplacian-Beltrami operator. 
Here we denoted by $\Lambda(\theta)$ the point of $\RR^{n+1}$ given 
by $(1,\theta)$. 
Indeed, we have
$$-\Delta_{\mathbb{H}^n}h_{\lambda,\theta}=
\left(\lambda^2+\frac{(n-1)^2}{4}\right)h_{\lambda,\theta}.$$
By analogy with the $\RR^n$ case, the definition of the Fourier transform is, 
for a function on $\mathbb{H}^n$,
$$\widehat{f}(\lambda,\theta)=
\int_{\mathbb{H}^n}h_{\lambda,\theta}(\Omega)
f(\Omega)\,d\Omega.$$
It turns out that this is the good definition, and one has the Fourier inversion formula for function on $\mathbb{H}^n$
$$f(\Omega)=\int_{-\infty}^\infty
\int_{\mathbb{S}^{n-1}}\overline{h}_{\lambda,\theta}(\Omega)
\widehat{f}(\lambda,\theta)\frac{d\theta d\lambda}{|c(\lambda)|^2},$$
where $c(\lambda)$ is the Harish-Chandra coefficient
$$\frac{1}{|c(\lambda)|^2}=\frac{1}{2(2\pi)^n}
\frac{|\Gamma(i\lambda+\frac{n-1}{2})|^2}{|\Gamma(i\lambda)|^2}.$$

Fore more details on hyperbolic space, see Helgason (\cite{Hel})
and Terras (\cite{Ter}).

\section{The representation of the solutions}

\subsection{The Euclidean case}

$$ $$

For the linear equation on $\RR^n$
$$\left\{ \begin{array}{c}
i\DP{t}u+\Delta_{\R^n} u=0,\\
u(0)=u_0\in L^2(\R^n),
\end{array}\right.,$$
the solution can be written explicitly. 
By applying the Fourier transform to the equation, one gets
$$\widehat{u}(\xi)=e^{-it|\xi|^2}\widehat{u_0}(\xi),$$
so the solution is
$$u(t,x)=
\frac{1}{(2\pi)^n}
\int_{\mathbb{R}^n}e^{-it|\xi|^2+ix\xi}\widehat{u_0}(\xi)d\xi.$$
By the Plancherel formula, the solution can also be written
$$u(t,x)=
\frac{1}{(4\pi it)^{\frac{n}{2}}}\int_{\mathbb{R}^n}
e^{i\frac{\mo{x-y}^2}{4t}}u_0(y)dy,$$
and therefore one obtains the following estimate
$$\|u(t)\|_\infty\leq \frac{c}{\mo{t}^\frac{n}{2}}\|u_0\|_1,$$
called dispersion estimate.

\subsection{Proof of Theorem \ref{exactsol}}
$$ $$
\par

For the hyperbolic space we will try the same approach, by using the Fourier transform introduced previously. By Fourier transforming on $\mathbb{H}^n$, one has

$$\left(i\partial_t-\left(\lambda^2+\frac{(n-1)^2}{4}\right)\right)
\widehat{u}(t,\lambda,\theta)=0,$$
so that
$$\widehat{u}(t,\lambda,\theta)=e^{-it\left(\lambda^2+\frac{(n-1)^2}{4}\right)}\widehat{u_0}(\lambda,\theta).$$
By applying the inverse Fourier transform,

$$u(t,\Omega)=\int_{-\infty}^\infty
\int_{\mathbb{S}^{n-1}}e^{-it(\lambda^2+\frac{(n-1)^2}{4})}
\overline{h}_{\lambda,\theta}(\Omega)
\widehat{u_0}(\lambda,\theta)\frac{d\theta d\lambda}{|c(\lambda)|^2}.$$
Finally, by making explicit the Fourier transform of the initial data,
\begin{equation}\label{kernel}
u(t,\Omega)=e^{-it\frac{(n-1)^2}{4}}\int_{\mathbb{H}^n}u_o(\Omega')\int_{-\infty}^\infty
e^{-it\lambda^2}L_\lambda(\Omega,\Omega') \frac{d\lambda}{|c(\lambda)|^2}d\Omega',\end{equation}
where
$$L_\lambda(\Omega,\Omega')=
\int_{\mathbb{S}^{n-1}}\overline{h}_{\lambda,\theta}(\Omega)h_{\lambda,\theta}(\Omega')d\theta.$$

\begin{lemma}\label{Lg}
For all $g\in SO(1,n)$, 
$$L_\lambda(\Omega,\Omega')=L_\lambda(g\Omega,g\Omega').$$
\end{lemma}

\begin{proof}

Using the definition (\ref{fp}) of the functions $h_{\lambda,\theta}$, 
and the fact that $g$ preserves the inner product, we have
\begin{equation}\label{identLg}
L_\lambda(\Omega,\Omega')=\int_{\mathbb{S}^{n-1}}[g\Omega,g\Lambda(\theta)]^{-i\lambda-\frac{n-1}{2}}
[g\Omega',g\Lambda(\theta)]^{i\lambda-\frac{n-1}{2}}d\theta.
\end{equation}

Consider now the map $F$ on $\mathbb{S}^{n-1}$,
$$F(\theta)=\tilde{\theta},$$
given by the relation
\begin{equation}\label{mu}
g\Lambda(\theta)=\mu(\theta)\Lambda(\tilde{\theta}),\end{equation}
for some real $\mu(\theta)$. 
Then, one can get (\cite{Vl}, Prop. 8.4.1) that this application is unique, is a 
diffeomorphism, 
$$\mu(\theta)=\cosh r_g+\sinh r_g \,\theta_1,$$
for some constant $r_g$ depending only on $g$, and the relation between the volume elements $d\tilde{\theta}$ and
$d\theta$ on $\mathbb{S}^{n-1}$ is
\begin{equation}\label{vol}
d\tilde{\theta}=\frac{d\theta}{(\cosh r_g+\sinh r_g\, \theta_1)^{n-1}}=\frac{d\theta}{\left(\mu(\theta)\right)^{n-1}}.\end{equation}
Here we denoted by $\theta_1$ the first $\mathbb{R}^n$-component of $\theta$.\\

Then, by making in the right hand side of (\ref{identLg}) 
the change of variable $F$, and by using (\ref{mu}) and
(\ref{vol}), we obtain 
$$L_\lambda(\Omega,\Omega')=\int_{\mathbb{S}^{n-1}}[g\Omega,\Lambda(\tilde{\theta})]^{-i\lambda-\frac{n-1}{2}}
[g\Omega',\Lambda(\tilde{\theta})]^{i\lambda-\frac{n-1}{2}}d\tilde{\theta}=
L_\lambda(g\Omega,g\Omega'),$$
so the Lemma is proved.
\end{proof}

Let us remark that Lemma \ref{Lg} can also be proved by using the classical
argument of the fundamental point pair invariant (\cite{D}, p.177).

\begin{lemma}\label{linkLegendre}
We have the following identity  
$$L_\lambda(\Omega,\Omega')=c (\sinh \rho)^{-\frac{n-2}{2}}
P^{-\frac{n-2}{2}}_{-\frac{1}{2}+i\lambda}(\cosh \rho),$$
where $\rho=d(\Omega,\Omega')$, and by $P^\mu_\nu(z)$ we denoted the Legendre
function, solution of the equation
$$(1-z^2)\partial^2_z w-2z\partial_z
w+\left(\nu(\nu+1)-\frac{\mu^2}{1-z^2}\right) w=0.$$
\end{lemma}

\begin{proof}

Let us choose a transformation $g\in SO(1,n)$, which maps $\Omega$ on the 
origin of $\mathbb{H}^n$. 
Then $g$ will send $\Omega'$ in a point of $\mathbb{H}^n$ of 
coordinates
$$\Omega^*=(\cosh \rho,\sinh \rho\,\gamma).$$
Since $g$ preserves the product, the radius of $\Omega^*$ can be calculated, 
$$\rho=\cosh^{-1}[\Omega,\Omega']=d(\Omega,\Omega').$$

So, by using Lemma \ref{Lg},
$$L_\lambda(\Omega,\Omega')=L_\lambda(0,\Omega^*)=
\int_{\mathbb{S}^{n-1}}[\Omega^*,\Lambda(\theta)]^{i\lambda-\frac{n-1}{2}}d\theta,$$
therefore
\begin{equation}\label{L}
L_\lambda(\Omega,\Omega')=
\int_{\mathbb{S}^{n-1}}(\cosh \rho-\sinh \rho\,\gamma.\theta)^{i\lambda-\frac{n-1}{2}}d\theta.
\end{equation}

Consider now a transformation $T\in SO(n)$ such that 
$$T(1,0,0...)=\gamma.$$
Then for $\theta\in\mathbb{S}^{n-1}$, we obtain an $\alpha\in[0,\pi[$ and a $\theta'\in\mathbb{S}^{n-2}$, well defined by the relation
$$\theta=T(\cos\alpha,\sin\alpha\,\theta').$$
Moreover, the angle between $\gamma$ and $\theta$ is $\alpha$,
$$\gamma.\theta=\cos \alpha.$$
The change of variable from $\theta$ to $\alpha$ and $\theta'$ is
$$d\theta=\sin^{n-2}\alpha \,d\alpha \,d\theta'.$$

Let us now apply this change of variables in (\ref{L}), and obtain
$$L_\lambda(\Omega,\Omega')=c\int_o^\pi (\cosh \rho-\sinh
\rho\,\cos\alpha)^{i\lambda-\frac{n-1}{2}}\sin^{n-2}\alpha \,d\alpha .$$
Therefore, in view of the integral form of $P^\mu_\nu(\cosh r)$
(\cite{B}), we obtain the identity of the Lemma.

\end{proof}

We shall now use the following Lemma, corresponding to Lemma 8.5.2 and
8.5.3 of \cite{Vl}.

\begin{lemma}\label{Legendre}
The function defined by 
$$L_\lambda^n(\rho)=c (\sinh \rho)^{-\frac{n-2}{2}}
P^{-\frac{n-2}{2}}_{-\frac{1}{2}+i\lambda}(\cosh \rho),$$
is, for $\rho>0$ and for $n\geq 1$ odd,  
\begin{equation}\label{Lodd}
L_\lambda^n(\rho)=c
\frac{|\Gamma(i\lambda)|^2}{|\Gamma(i\lambda+\frac{n-1}{2})|^2}\left(\frac{\partial_\rho}{\sinh\rho}\right)^\frac{n-1}{2}\cos\lambda\rho,
\end{equation}
and for $n\geq 2$ even,  
\begin{equation}\label{Leven}
L_\lambda^n(\rho)=c
\frac{|\Gamma(i\lambda)|^2}{|\Gamma(i\lambda+\frac{n-1}{2})|^2}
\int_\rho^\infty\frac{\sinh s}{\sqrt{\cosh s-\cosh\rho}}
\left(\frac{\partial_s}{\sinh s}\right)^\frac{n}{2}\cos\lambda s\, ds.
\end{equation}
\end{lemma}

By using in the formula (\ref{kernel}) the expression of
$L_\lambda(\Omega,\Omega')$ given by the Lemma \ref{linkLegendre} and
Lemma \ref{Legendre}, we obtain that
$$u(t,\Omega)=ce^{-it\frac{(n-1)^2}{4}}\int_{\mathbb{H}^n}u_o(\Omega')K^n(t,d(\Omega,\Omega'))d\Omega',$$
where the kernel $K^n$ is, for$\rho>0$ and for $n\geq 3$ odd,

$$K^n(t,\rho)=\int_{-\infty}^\infty
e^{-it\lambda^2}
\left(\frac{\partial_\rho}{\sinh\rho}\right)^\frac{n-1}{2}\cos\lambda\rho\,
d\lambda,$$
and for $n\geq 2$ even,  
\begin{equation}\label{Keven}
K^n(t,\rho)=\int_{-\infty}^\infty e^{-it\lambda^2}
\int_\rho^\infty\frac{\sinh s}{\sqrt{\cosh s-\cosh\rho}}
\left(\frac{\partial_s}{\sinh s}\right)^\frac{n}{2}\cos\lambda s\, ds\,
d\lambda.
\end{equation}

\par
In the case of odd dimensions, the kernel can be written
\begin{equation}\label{Kodd}
K^n(t,\rho)=\left(\frac{\partial_\rho}{\sinh\rho}\right)^\frac{n-1}{2}
\int_{-\infty}^\infty
e^{-it\lambda^2}\cos\lambda\rho\,
d\lambda=\frac{c}{|t|^\frac{1}{2}}\left(\frac{\partial_\rho}{\sinh\rho}\right)^\frac{n-1}{2}e^{i\frac{\rho^2}{4t}},\end{equation}
so the representation (\ref{odd}) of the Theorem \ref{exactsol} is obtained.

\par
In the case of even dimensions, the kernel can be expressed by
induction (\cite{Vl}, formula (8.5.23)) in terms of the kernel 
of the dimension $2$
$$K^n(t,\rho)=\left(\frac{\partial_\rho}{\sinh
\rho}\right)^\frac{n-2}{2}K^2(t,\rho).$$
As $K^2$ can be written, after applying the Fubini theorem in
(\ref{Keven}), in the form
$$K^2(t,\rho)=\int_\rho^\infty\frac{1}{\sqrt{\cosh s-\cosh\rho}}
\int_{-\infty}^\infty e^{-it\lambda^2}\sin \lambda s\,\lambda \,d\lambda\,ds=
\frac{c}{|t|^\frac{3}{2}}\int_\rho^\infty\frac{e^{i\frac{s^2}{4t}}s}
{\sqrt{\cosh s-\cosh\rho}}ds,$$
we get that
\begin{equation}\label{kerpair}
K^n(t,\rho)=\frac{c}{|t|^\frac{3}{2}}\left(\frac{\partial_\rho}{\sinh
\rho}\right)^\frac{n-2}{2}\int_\rho^\infty\frac{e^{i\frac{s^2}{4t}}s}{\sqrt{\cosh
s-\cosh\rho}}ds,\end{equation}
so we obtain also the representation (\ref{even}) and the first point
of Theorem \ref{exactsol} is proved.

\section{The $L^1-L^\infty$ dispersion estimates}

We shall use the following Proposition, which shall be proved in the
Appendix A.
\begin{prop}\label{F}
For all $m\geq 1$ integer, the following identity holds 
\begin{equation}\label{derexp}
\left(\frac{\partial_s}{\sinh s}\right)^me^{i\frac{s^2}{4t}}=
\sum_{k=1}^m \frac{e^{i\frac{s^2}{4t}}}{t^k}\,F_k^m(s),\end{equation}
where, modulo constants,  
\begin{equation}\label{idFm}
F_m^m(s)=\left(\frac{s}{\sinh s}\right)^m,\end{equation}
and, for $1\leq k<m$, $F_k^m(s)$ can be written
\begin{equation}\label{idF}
F_k^m(s)=\sum_{r=0}^k\,\,\sum_{\left\{\begin{array}{c}i_1+\cdot\cdot\cdot+i_r=m-k\\i_j\geq 1\end{array}\right.}\left(\frac{s}{\sinh s}\right)^{k-r}
\left(\frac{\partial_s}{\sinh s}\right)^{i_1}
\frac{s}{\sinh s}...
\left(\frac{\partial_s}{\sinh s}\right)^{i_r}
\frac{s}{\sinh s}.\end{equation}
Moreover, for $\alpha\in\{0,1\}$ and $1\leq k\leq m$
\begin{equation}\label{estFk}
\left|\partial_s^\alpha\left(\frac{\sinh s}{s}\,F_k^m(s)\right)\right|\leq 
c\,s^\alpha\left(\frac{s}{\sinh s}\right)^{m-1},\end{equation}
and the second derivative of $\frac{\sinh s}{s}\,F_k^m(s)$ is bounded for $m\geq 2$.

\end{prop}
\subsection{Dispersive estimates for odd dimensions} 

From the expression (\ref{Kodd}), for odd dimensions, the kernel is for $\rho>0$
$$K^n(t,\rho)=\frac{c}{|t|^\frac{1}{2}}\left(\frac{\partial_\rho}{\sinh\rho}\right)^\frac{n-1}{2}e^{i\frac{\rho^2}{4t}},$$
so by (\ref{derexp}) it can be developed in powers of $t$ as follows
\begin{equation}\label{kerneldevp}
K^n(t,\rho)=\sum_{k=1}^{\frac{n-1}{2}} \frac{e^{i\frac{s^2}{4t}}}{t^{k+\frac{1}{2}}}F_k^{\frac{n-1}{2}}(\rho).
\end{equation}
Combining this with (\ref{estFk}) for $\alpha =0$ we obtain the dispersion estimates
(\ref{disph}) and (\ref{disphtlarge}) for all $n\geq 3$ odd.\\

\subsection{Local dispersion for even dimensions} 
Between two kernels of consecutive order, there is a induction
relation, namely
$$K^n(t,\rho)=c\int_\rho^\infty \frac{\sinh s}{\sqrt{\cosh
s-\cosh\rho}}\,K^{n+1}(t,s)\,ds,$$
(see for example formula ($5.7.4$) in the book of Davies (\cite{D})).

Therefore, by using the development (\ref{kerneldevp}) of the kernels
of odd order, we can get the one for kernels of even order for $\rho>0$
\begin{equation}\label{kernel2devp} 
K^n(t,\rho)=\sum_{k=1}^{\frac{n}{2}}\frac{1}{t^{k+\frac{1}{2}}}\int_\rho^\infty 
\frac{e^{i\frac{s^2}{4t}}s}{\sqrt{\cosh s-\cosh \rho}}\frac{\sinh s}{s} F_k^\frac{n}{2}(s)\,ds.\end{equation}
Therefore the local dispersion (\ref{disph}) is obtained for $n=2$ by (\ref{idFm}) and the following Proposition. For all even dimensions $n\geq 4$, the estimates (\ref{mainoscintegral}) and (\ref{oscintegral}) of Corollary \ref{corol} give us
$$|K^n(t,\rho)|\leq \left(\frac{\rho}{\sinh \rho}\right)^\frac{n-1}{2}\left(\frac{1}{t^{\frac{n}{2}}}+\sum_{k=1}^{\frac{n-2}{2}}\frac{1}{t^{k+\frac{1}{2}}}\right),$$
and the local dispersion (\ref{disph}) follows.

\begin{prop}\label{decay}
The integral
$$I(t,\rho)=\int_\rho^\infty\frac{e^{i\frac{s^2}{4t}}s}{\sqrt{\cosh s-\cosh\rho}}ds$$
can be estimated, for small times, by
\begin{equation}\label{decayappendix}
|I(t,\rho)|\leq c\sqrt{t}\sqrt{\frac{\rho}{\sinh \rho}}.
\end{equation}
\end{prop}

\begin{proof} 
In the following we shall use frequently the estimate
\begin{equation}\label{devlimite}
\frac{1}{\sqrt{\cosh s-\cosh\rho}}\leq
\frac{c}{(s-\rho)\sqrt{\cosh\rho}}\leq \frac{c}{s-\rho}
\end{equation}
for $s>\rho\geq 0$, and
\begin{equation}\label{devlimite2}
\frac{1}{\sqrt{\cosh s-\cosh\rho}}\leq \frac{c}{\sqrt{(s-\rho)\sinh\rho}}
\end{equation}
for $s>\rho>0$.\\

For small time, in order to get the decay (\ref{decayappendix}), 
we have to get advantage of the imaginary phase. Let $t<\frac{1}{2}$.

{\underline{The case $\rho\geq\frac{\sqrt{t}}{2}$}}.\\

By doing in the integral $I(t,\rho)$ the change of variable
$s=\tau\frac{t}{\rho}+\rho$, we obtain
$$I(t,\rho)=\frac{te^{i\frac{\rho^2}{4t}}}{\rho}
\int_0^\infty
\frac{e^{i\tau^2\frac{t}{4\rho^2}+i\frac{\tau}{2}}(\tau\frac{t}{\rho}+\rho)}
{\sqrt{\cosh (\tau\frac{t}{\rho}+\rho)-\cosh\rho}}d\tau=$$
$$=2e^{i\frac{\rho^2}{4t}}\sqrt{t}\sqrt{\frac{\rho}{\sinh\rho}}
\int_0^\infty
\frac{e^{i\tau^2\frac{t}{4\rho^2}+i\frac{\tau}{2}}(\tau\frac{t}{2\rho^2}+\frac{1}{2})}
{\sqrt{\frac{\rho}{t\sinh\rho}\left(\cosh (\tau\frac{t}{\rho}+\rho)-\cosh\rho\right)}}d\tau.$$
We split the integral into two pieces
\begin{equation}\label{unifbound}
I(t,\rho)=2e^{i\frac{\rho^2}{4t}}\sqrt{t}\sqrt{\frac{\rho}{\sinh\rho}}
\left(\int_0^1 \, +\int_1^\infty\right)=2e^{i\frac{\rho^2}{4t}}\sqrt{t}\sqrt{\frac{\rho}{\sinh\rho}}
(J_1+J_2).
\end{equation}
We shall show that $J_1$ and $J_2$ are bounded independently of $t$ and
$\rho$.

By using (\ref{devlimite2}), 
$$|J_1|\leq\frac{1}{2}\int_0^1\frac{\tau\frac{t}{\rho^2}+1}
{\sqrt{\frac{\rho}{t\sinh\rho}\tau\frac{t}{\rho}\sinh\rho}}=
\frac{1}{2}\int_0^1\frac{\tau\frac{t}{\rho^2}+1}{\sqrt{\tau}}.$$
Since we are in the case $\rho\geq\frac{\sqrt{t}}{2}$, it follows that
$$\frac{t}{\rho^2}\leq 4,$$
and therefore 
$$|J_1|\leq\frac{1}{2}\int_0^1\frac{5}{\sqrt{\tau}}=5,$$
The second term can be written
$$J_2=\int_1^\infty e^{i\phi (\tau)}\phi'(\tau)\alpha(\tau)d\tau,$$
where
$$\phi(\tau)=\tau^2\frac{t}{4\rho^2}+\frac{\tau}{2},$$
and 
$$\alpha(\tau)=\frac{1}{\sqrt{\frac{\rho}{t\sinh\rho}\left(\cosh (\tau\frac{t}{\rho}+\rho)-\cosh\rho\right)}}.$$
By integrating by parts in $J_2$ we get
$$J_2=\frac{e^{i\phi (\tau)}}{i}\alpha(\tau)|_1^\infty-
\int_1^\infty \frac{e^{i\phi (\tau)}}{i}\alpha'(\tau)d\tau,$$
so
\begin{equation}\label{I0}
|J_2|\leq 2\,\underset{1\leq\tau}{\sup}|\alpha|+\int_1^\infty |\alpha'(\tau)|d\tau.\end{equation}

Let us notice that by (\ref{devlimite2})
$$\alpha(1)\leq 
\frac{1}{\sqrt{\frac{\rho}{t\sinh\rho}\frac{t}{\rho}\sinh\rho}}=1,$$
and that $\alpha$ is a decreasing function. Therefore,
$$|J_2|\leq 2+\int_1^\infty |\alpha'(\tau)|d\tau.$$
Since the derivative of $\alpha$ is negative,
$$|J_2|\leq 2+\left|\int_1^\infty \alpha'(\tau)d\tau\right|
=2+\left| \alpha(\tau)|_1^\infty\right|\leq 4.$$

Therefore we have obtained that $J_1$ and $J_2$ are bounded independently of $t$ and
$\rho$. In view of (\ref{unifbound}) the decay (\ref{decayappendix}) is obtained in the region 
$\rho\geq\frac{\sqrt{t}}{2}$.

{\underline{The case $\rho<\frac{\sqrt{t}}{2}$}}.\\

Since $t<\frac{1}{2}$, then $\rho<\frac{1}{2\sqrt{2}}$, and by noticing that the quotient $\frac{\rho}{\sinh\rho}$ is bounded 
near zero, it will be sufficient to prove that
\begin{equation}\label{small} 
\int_\rho^\infty\frac{e^{i\frac{s^2}{4t}}s}{\sqrt{\cosh
s-\cosh\rho}}ds\leq c\sqrt{t}.
\end{equation}
Let us split this integral in three parts :
$$\int_\rho^\infty\frac{e^{i\frac{s^2}{4t}}s}{\sqrt{\cosh s-\cosh\rho}}ds=
\int_\rho^{2\rho}+\int_{2\rho}^{\sqrt{t}}+\int_{\sqrt{t}}^\infty=
I_1+I_2+I_3.$$

If $\rho=0$ then $I_1$ is zero, otherwise by (\ref{devlimite2})
$$|I_1|\leq \int_\rho^{2\rho}\frac{s}{\sqrt{\cosh s-\cosh\rho}}ds\leq 
2\rho \int_\rho^{2\rho}\frac{1}{\sqrt{(s-\rho)\sinh \rho}}ds= $$
$$=4\rho\frac{\sqrt{\rho}}{\sqrt{\sinh \rho}}.$$
Since we are in the case $\rho<\frac{\sqrt{t}}{2}$ we get that
$$|I_1|\leq 2\sqrt{t}\sqrt{\frac{\rho}{\sinh \rho}}\leq c\sqrt{t}.$$

If $s\geq 2\rho$, which is the case in $I_2$, then by (\ref{devlimite})
$$\frac{s}{\sqrt{\cosh s-\cosh\rho}}\leq \frac{s}{s-\rho}\leq 3,$$
and we obtain
$$|I_2|\leq 3(\sqrt{t}-2\rho)\leq 3\sqrt{t}.$$

By performing in $I_3$ the change of variable
$s=\sqrt{t}\tau$,
$$I_3=\sqrt{t}\int_1^\infty e^{i\frac{\tau^2}{4}}\,a(\sqrt{t}\tau)d\tau,$$
where 
$$a(s)=\frac{s}{\sqrt{\cosh s-\cosh\rho}}$$
for $s\geq\sqrt{t}$. By integrating by parts two times the integral $I_3$, we get
$$|I_3|\leq c\sqrt{t},$$
provided that the following lemma holds. (see
Theorem I.8.1 of \cite{AG}).

\begin{lemma}\label{ab}
The function $a(\sqrt{t}\tau)$, its first and second
derivative in $\tau$, are all bounded independently of $t$, for all 
$\tau\geq 1$.
\end{lemma}

\begin{proof}
By using the same argument as in the estimate of $I_2$, we have that
$a(\sqrt{t}\tau)$ is upper-bounded by $3$ for $\tau\geq 1$.

The derivative of $a$ is
$$a'(s)=\frac{1}{\sqrt{\cosh s-\cosh\rho}}-
\frac{s\sinh s}{2(\cosh s-\cosh\rho)^{\frac{3}{2}}}.$$
Since we are in the region $s\geq\sqrt{t}$,
$\rho<\frac{\sqrt{t}}{2}$, the first term can be estimated by (\ref{devlimite})
$$\frac{1}{\sqrt{\cosh s-\cosh\rho}}\leq \frac{c}{s-\rho}\leq \frac{c}{\sqrt{t}-\frac{\sqrt{t}}{2}}
\leq \frac{c}{\sqrt{t}}.$$
The second term in the derivative of $a$ 
can be written as 
$$-\frac{a(s)}{2}\frac{\sinh s}{\cosh s-\cosh\rho}.$$
Since $s\geq\sqrt{t}>\rho$, the second fraction is positive, and as $a$ was already proved to be bounded, we get 
$$|a'(s)|\leq \frac{c}{\sqrt{t}}+c\,\underset{s\geq\sqrt{t}}{\sup}\,\,\frac{\sinh s}{\cosh s-\cosh\rho}.$$
We have
$$\left(\frac{\sinh s}{\cosh s-\cosh\rho}\right)'=\frac{1-\cosh s\cosh\rho}{\left(\cosh s-\cosh\rho\right)^2}\leq 0,$$
so it follows that
$$\underset{s\geq\sqrt{t}}{\sup}\,\,\frac{\sinh s}{\cosh s-\cosh\rho}=\frac{\sinh \sqrt{t}}{\cosh \sqrt{t}-\cosh\rho}.$$
As $t<\frac{1}{2}$ and $\rho<\frac{\sqrt{t}}{2}$,
\begin{equation}\label{decr}
\underset{s\geq\sqrt{t}}{\sup}\,\,\frac{\sinh s}{\cosh s-\cosh\rho}\leq c\frac{\sqrt{t}}{(\sqrt{t}-\rho)^2}\leq c\frac{\sqrt{t}}{(\sqrt{t}-\frac{\sqrt{t}}{2})^2}\leq \frac{c}{\sqrt{t}},\end{equation}
and we have obtained that for $s\geq\sqrt{t}$
$$|a'(s)|\leq \frac{c}{\sqrt{t}}.$$
It follows that
$$\partial_\tau a(\sqrt{t}\tau)=\sqrt{t}a'(\sqrt{t}\tau)$$
is bounded independently of $t$ for $\tau\geq 1$.

Finally, the second derivative of $a$ is
$$a''(s)=-\frac{\sinh s}{2(\cosh s-\cosh\rho)^{\frac{3}{2}}}-\frac{s\cosh s}{(\cosh s-\cosh\rho)^{\frac{3}{2}}}+
\frac{3s\sinh^2s}{4(\cosh s-\cosh\rho)^{\frac{5}{2}}}.$$
By using (\ref{decr}) and then (\ref{devlimite}) we can treat the first term
$$\frac{\sinh s}{2(\cosh s-\cosh\rho)^{\frac{3}{2}}}\leq \frac{c}{\sqrt{t}}\frac{1}{\sqrt{\cosh s-\cosh\rho}}\leq \frac{c}{\sqrt{t}}\frac{c}{s-\rho}\leq \frac{c}{t}.$$
The second term in the derivative can be written
$$\frac{s\cosh s}{(\cosh s-\cosh\rho)^{\frac{3}{2}}}=a(s)\frac{\cosh s}{\cosh s-\cosh\rho}=a(s)\left(1+\frac{\cosh \rho}{\cosh s-\cosh\rho}\right).$$
As $a$ was proved to be bounded and as the fraction in the brackets is decreasing, for $s\geq \sqrt{t}$
$$\frac{s\cosh s}{(\cosh s-\cosh\rho)^{\frac{3}{2}}}\leq c\left(1+\frac{\cosh\rho}{(\sqrt{t}-\rho)^2}\right).$$
We are in the region $\rho<\frac{\sqrt{t}}{2}<\frac{1}{2\sqrt{2}}$, so the second term in the derivative can be estimated
$$\frac{s\cosh s}{(\cosh s-\cosh\rho)^{\frac{3}{2}}}\leq c\left(1+\frac{c}{t}\right)\leq \frac{c}{t}.$$
Finally, the third term in the derivative can be written as
$$\frac{3a(s)}{4}\left(\frac{\sinh s}{\cosh s-\cosh\rho}\right)^2,$$
and can be upper-bounded using (\ref{decr}) by $\frac{c}{t}$.
In conclusion, for $s\geq\sqrt{t}$
$$|a''(s)|\leq \frac{c}{t}.$$
So we get that 
$$\partial_\tau ^2a(\sqrt{t}\tau)=ta''(\sqrt{t}\tau)$$
is bounded independently of $t$ for $\tau\geq 1$.
\end{proof}

In conclusion, we have obtained that the three integrals $I_1, I_2$
and $I_3$ are upper-bounded by $c\sqrt{t}$, and (\ref{small})
follows. Therefore the estimate (\ref{decayappendix}) is proved and 
Proposition \ref{decay} follows. 

\end{proof}

\begin{corollary}\label{corol}
For small times and $n\geq 4$ the following estimates hold
\begin{equation}\label{mainoscintegral}
\left|\int_\rho^\infty\frac{e^{i\frac{s^2}{4t}}s}{\sqrt{\cosh s-\cosh\rho}}
\frac{\sinh s}{s}F^{\frac{n}{2}}_{\frac{n}{2}}(s)\,ds\right|\leq 
c\sqrt{t}\left(\frac{\rho}{\sinh \rho}\right)^{\frac{n-1}{2}},
\end{equation}
and for $1\leq k< \frac{n}{2}$
\begin{equation}\label{oscintegral}
\left|\int_\rho^\infty\frac{e^{i\frac{s^2}{4t}}s}{\sqrt{\cosh s-\cosh\rho}}
\frac{\sinh s}{s}F^{\frac{n}{2}}_k(s)\,ds\right|\leq 
c\left(\frac{\rho}{\sinh \rho}\right)^{\frac{n-1}{2}}.\end{equation}
\end{corollary}

\begin{proof}
We redo the proof of Proposition \ref{decay}, which is the particular
case $n=2$ of its Corollary. The only delicate point will be the integral $J_2$ which will generate 
the two different estimates (\ref{mainoscintegral}) and (\ref{oscintegral}).

By using the upper-bound (\ref{estFk}) for $\alpha =0$ and $1\leq k\leq\frac{n}{2}$, and the fact that the integration is done for $s\geq\rho$,
$$\left|\frac{\sinh s}{s}F^{\frac{n}{2}}_k(s)\right|\leq \left(\frac{s}{\sinh s}\right)^{\frac{n-2}{2}}\leq \left(\frac{\rho}{\sinh \rho}\right)^{\frac{n-2}{2}}.$$
Therefore everytime the imaginary phase is ignored, we find ourselves in the same situation as in Proposition \ref{decay}, and we obtain the desired estimates.

The imaginary phase has been taken in account only in the terms $J_2$
and $I_3$.

For estimating $I_3$ we have to prove Lemma \ref{ab} with
$$\widetilde{a}(s)=a(s)\frac{\sinh s}{s}F^{\frac{n}{2}}_k(s).$$
We had that $a$ is bounded by a constant, its first derivative bounded
by $ct^{-\frac{1}{2}}$, and its second derivative bounded by
$ct^{-1}$. We want the same for $\widetilde{a}$. As $t$ is small, it will be sufficient to prove that $\frac{\sinh s}{s}F^{\frac{n}{2}}_k(s)$ is
bounded, and also are its first and second derivative. Since we are in the case $n\geq 4$, Proposition \ref{F} tells us already that its second derivative is bounded. Moreover, for $1\leq k\leq\frac{n}{2}$, by (\ref{estFk}) for $\alpha =0$ 
$$\left|\frac{\sinh s}{s}F^{\frac{n}{2}}_k(s)\right|\leq \left(\frac{s}{\sinh
s}\right)^{\frac{n-2}{2}},$$
and by (\ref{estFk}) for $\alpha =1$ 
$$\left|\partial_s\frac{\sinh s}{s}F^{\frac{n}{2}}_k(s)\right|\leq s\left(\frac{s}{\sinh
s}\right)^{\frac{n-2}{2}}=\frac{s^2}{\sinh s}\left(\frac{s}{\sinh
s}\right)^{\frac{n-4}{2}}.$$
Again, since we are in the case $n\geq 4$, these two quantities are bounded, so Lemma \ref{ab} still holds in the context of the Corollary.

In conclusion, up to the term $J_2$, our integrals have the decay of Proposition \ref{decay}
$$\sqrt{t}\sqrt{\frac{\rho}{\sinh \rho}}$$
improved by the extra term
$$\left(\frac{\rho}{\sinh \rho}\right)^{\frac{n-2}{2}}.$$

Let us now first treat $J_2$ for obtaining (\ref{mainoscintegral}). By (\ref{idFm}), 
$$\frac{\sinh s}{s}F^{\frac{n}{2}}_{\frac{n}{2}}(s)=\left(\frac{s}{\sinh
s}\right)^{\frac{n-2}{2}}.$$
Therefore the same arguments as in Proposition \ref{decay} can be performed in $J_2$, by replacing the
function $\alpha(\tau)$ by

$$\alpha(\tau)\left(\frac{\tau\frac{t}{\rho}+\rho}{\sinh\left(\tau\frac{t}{\rho}+\rho\right)}\right)^{\frac{n-2}{2}},$$
which is still a decreasing function, bounded at $1$ by 
$\left(\frac{\rho}{\sinh\rho}\right)^{\frac{n-2}{2}}$. So $J_2$ will
be also upper-bounded by $\left(\frac{\rho}{\sinh\rho}\right)^{\frac{n-2}{2}}$ and (\ref{mainoscintegral}) is proved.

In view to obtain (\ref{oscintegral}), we restart the argument performed on $J_2$ in Proposition \ref{decay}, with $\alpha(\tau)$ replaced by
$$\alpha(\tau)\frac{\sinh s}{s}F^{\frac{n}{2}}_k(s),$$
where
$$s=\tau\frac{t}{\rho}+\rho.$$
When getting to (\ref{I0}) we need to estimate
$$I_4=\int_1^\infty\left|\partial_\tau\left(\alpha(\tau)\frac{\sinh s}{s}F^{\frac{n}{2}}_k(s)\right)\right|d\tau.$$
Now the function $\frac{\sinh s}{s}F^{\frac{n}{2}}_k(s)$ is not necessarily decreasing, so we are not able to get rid of the modulus in the integral, as easily as before. We will be able to estimate  
\begin{equation}\label{I4}
I_4\leq\frac{1}{\sqrt{t}}\left(\frac{\rho}{\sinh\rho}\right)^{\frac{n-2}{2}},\end{equation}
that is with a loss of $\sqrt{t}$.\\
In the derivation of the product,
$$I_4=
\int_1^\infty\left|\alpha'(\tau)\frac{\sinh s}{s}F^{\frac{n}{2}}_k(s)+\alpha(\tau)\frac{t}{\rho}\partial_s\left(\frac{\sinh s}{s}F^{\frac{n}{2}}_k(s)\right)\right|d\tau,$$
we can use the estimate (\ref{estFk}) for $\alpha\in\{0,1\}$ and then
$$I_4\leq
\left(\frac{\rho}{\sinh\rho}\right)^{\frac{n-2}{2}}\int_1^\infty\left|\alpha'(\tau)\right|d\tau+
\left(\frac{\rho}{\sinh\rho}\right)^{\frac{n-2}{2}}\int_1^\infty\left|
\alpha(\tau)\frac{t}{\rho}\left(\tau\frac{t}{\rho}+\rho\right)\right|d\tau.$$
The first integral is exactly the one in Proposition \ref{decay}, so it has been proved that it is bounded independently of $t$ and $\rho$. For treating the second integral, we shall use that for all odd $N$,
$$\alpha(\tau)=\frac{1}{\sqrt{\cosh \left(\tau\frac{t}{\rho}+\rho\right)-\cosh \rho}}\leq 
\frac{1}{\sqrt{\sinh\rho}}\left(\tau\frac{t}{\rho}\right)^{-\frac{N}{2}}.$$
Using this estimate for $N=5$ and for $N=3$,
$$\int_1^\infty\left|
\alpha(\tau)\frac{t}{\rho}\left(\tau\frac{t}{\rho}+\rho\right)\right|d\tau\leq 
\frac{1}{\sqrt{\sinh\rho}}\int_1^\infty \left(\left(\tau\frac{t}{\rho}\right)^{-\frac{5}{2}}\tau\left(\frac{t}{\rho}\right)^2+\left(\tau\frac{t}{\rho}\right)^{-\frac{3}{2}}\frac{t}{\rho}\rho\right)d\tau.$$
The integrals in $\tau$ are bounded, so
$$\int_1^\infty\left|\alpha(\tau)\frac{t}{\rho}\left(\tau\frac{t}{\rho}+\rho\right)\right|d\tau\leq 
\frac{1}{\sqrt{\sinh\rho}}\left(\frac{t}{\rho}\right)^{-\frac{1}{2}}(1+\rho)\leq \frac{c}{\sqrt{t}},$$
so estimate (\ref{I4}) follows and the (\ref{oscintegral}) is proved.
\end{proof}

\subsection{Large time dispersion for even dimensions}

By using the development (\ref{kernel2devp}) of the dispersion kernel and the estimates (\ref{estFk}) for $\alpha =0$ we obtain
\begin{equation}\label{kernel2devplarget} 
|K^n(t,\rho)|\leq \sum_{k=1}^\frac{n}{2}\frac{1}{t^{k+\frac{1}{2}}}\left(\frac{\rho}{\sinh \rho}\right)^\frac{n-2}{2}\left|\int_\rho^\infty 
\frac{s}{\sqrt{\cosh s-\cosh \rho}}\,ds\right|.\end{equation}
Let us notice that we have not taken into account the oscillatory phase in 
the integral. We shall split the
remaining integral into two pieces
$$\left| \int_\rho^\infty\frac{s}{\sqrt{\cosh s-\cosh\rho}}ds \right|=
\left| \int_0^\infty\frac{s+\rho}{\sqrt{\cosh (s+\rho)-\cosh\rho}}ds \right|=
\int_0^1+\int_1^\infty.$$
For $\rho>0$, the first part is upper-bounded by
$$\int_0^1\frac{s+\rho}{\sqrt{s\sinh\rho}}ds\leq c\frac{1+\rho}{\sqrt{\sinh\rho}},$$
and the second, for a $N$ odd integer sufficiently large, by
$$\int_1^\infty\frac{s+\rho}{\sqrt{s^N\sinh\rho}}ds\leq c\frac{1+\rho}{\sqrt{\sinh\rho}}.$$
Now, by (\ref{kernel2devplarget}) we get for $\rho>0$
$$|K^n(t,\rho)|\leq \sum_{k=1}^\frac{n}{2}\frac{1}{t^{k+\frac{1}{2}}}\left(\frac{\rho}{\sinh \rho}\right)^\frac{n-2}{2}\frac{1+\rho}{\sqrt{\sinh\rho}}.$$
Therefore we have obtained the dispersion estimate (\ref{disphtlarge2}) for large times, and the last part 
of Theorem \ref{disp} i) follows.

\section{The weighted dispersion estimates}
Let us prove now the weighted dispersion (\ref{dispw}). From
(\ref{disph}) we have that for small times and $n\geq 3$
\begin{equation}\label{ppp}
|u(t,\Omega)|\leq \frac{c}{|t|^\frac{n}{2}}
\int_{\mathbb{H}^n}|u_0(\Omega')|\left(\frac{\rho}{\sinh\rho}\right)^\frac{n-1}{2}\,d\Omega'\leq \frac{c}{|t|^\frac{n}{2}}
\int_{\mathbb{H}^n}|u_0(\Omega')|\frac{\rho}{\sinh\rho}\,d\Omega'.
\end{equation}
The initial data is considered here radial, that is
$$u_0(\Omega')=u_0(\cosh r', \sinh r'\omega')=u_0(r').$$
Also in hyperbolic coordinates we can write, using (\ref{prod}) and
(\ref{dist})
$$\rho=d(\Omega,\Omega')=\cosh^{-1}[\Omega,\Omega']=
\cosh ^{-1}(\cosh r\cosh r'-\sinh r\sinh r'\,\omega.\omega').$$
Therefore, by passing in hyperbolic coordinates in (\ref{ppp}),
\begin{equation}\label{qqq}
|u(t,\Omega)|=|u(t,r)|\leq \frac{c}{|t|^\frac{n}{2}}
\int_0^\infty |u_0(r')|K(t,r,r')\sinh^{n-1}r'\,dr',\end{equation}
where
$$K(t,r,r')=\int_{\mathbb{S}^{n-1}}
\frac{\cosh ^{-1}(\cosh r\cosh r'-\sinh r\sinh r'\,\omega.\omega')}
{\sinh \cosh ^{-1}(\cosh r\cosh r'-\sinh r\sinh
r'\,\omega.\omega')}\,d\omega'.$$
By doing a rotation as in the proof of the Lemma \ref{linkLegendre},
$$K(t,r,r')=c\int_0^\pi 
\frac{\cosh ^{-1}(\cosh r\cosh r'-\sinh r\sinh r'\,\cos\alpha)}
{\sinh \cosh ^{-1}(\cosh r\cosh r'-\sinh r\sinh
r'\,\cos\alpha)}\sin^{n-2}\alpha\,d\alpha.$$
Now let us do the change of variable $\cos\alpha=x$ and get
$$K(t,r,r')=c\int_{-1}^1 
\frac{\cosh ^{-1}(\cosh r\cosh r'-\sinh r\sinh r'\,x)}
{\sinh \cosh ^{-1}(\cosh r\cosh r'-\sinh r\sinh
r'\,x)}(1-x^2)^\frac{n-3}{2}dx.$$
Since here $n$ is larger or equal to $3$,
$$K(t,r,r')\leq c\int_{-1}^1 
\frac{\cosh ^{-1}(\cosh r\cosh r'-\sinh r\sinh r'\,x)}
{\sinh \cosh ^{-1}(\cosh r\cosh r'-\sinh r\sinh
r'\,x)}dx.$$
Finally, let us do the change of variable
$$\cosh r\cosh r'-\sinh r\sinh r'\,x=\cosh y,$$
and then
$$K(t,r,r')\leq \frac{c}{\sinh r\sinh r'}\int_{|r-r'|}^{r+r'} y\,dy=
c\frac{r\,r'}{\sinh r\sinh r'}.$$
Using now (\ref{qqq}), we obtain that
$$|u(t,\Omega)|=|u(t,r)|\leq \frac{c}{|t|^\frac{n}{2}}
\int_0^\infty |u_0(r')|\frac{r\,r'}{\sinh r\sinh r'}\sinh^{n-1}r'\,dr'.$$
Remembering that $r=d(\Omega,0)$ we can go back to the
$\Omega$-coordinates and
$$|u(t,\Omega)|\,\frac{\sinh d(\Omega,0)}{d(\Omega,0)}\leq \frac{c}{|t|^\frac{n}{2}}
\int_{\mathbb{H}^n}|u_0(\Omega')|\,\frac{d(\Omega',0)}{\sinh
d(\Omega',0)}d\Omega',$$
so the weighted local dispersion estimate (\ref{dispw}) is proved.\\

For proving the weighted Strichartz estimates we will apply the general
lemma of Keel and Tao (Theorem 10.1 of \cite{KT}), in the case 
$$H=B_0=L^2_{rad}\;\;,\;\;B_1=L^1_{rad}(w^{-1})\;\;,\;\;
\sigma=\frac{n}{2}\;\;,\;\;2\leq \frac{4}{n\theta}\;\;,\;\;
0\leq\theta\leq 1 \;\;,\;\;\left(\frac{4}{n\theta},\theta,\sigma\right)\neq (2,1,1).$$
We obtain that for radial $u$
$$\|u\|_{L^\frac{4}{n\theta}\left([0,T],\,
(L^2,L^1(w^{-1}))^*_{\theta,2}\right)}\leq 
c\|u_0\|_{L^2},$$
where $(L^2,L^1(w^{-1}))^*_{\theta,2}$ is the dual of the real
interpolation space between $L^2$ and $L^1(w^{-1})$. The conditions on $\theta$ can be rewritten, for $n\geq 2$,
$$0\leq \theta\leq \frac{2}{n}\;\;,\;\;(\theta,n)\neq (1,2).$$

For any $r\leq 2$ we have
$$(L^2,L^1(w^{-1}))_{\theta,r}\subseteq (L^2,L^1(w^{-1}))_{\theta,2},$$
so that
\begin{equation}\label{rrr}
\|u\|_{L^\frac{4}{n\theta}\left([0,T],\,
(L^2,L^1(w^{-1}))^*_{\theta,r}\right)}\leq 
c\|u_0\|_{L^2}.\end{equation}
By the Theorem 5.5.1. of \cite{BL}, for $0<\theta<1$, 
$$(L^2,L^1(w^{-1}))_{\theta,r}=L^r(\widetilde{w}),$$
with
$$\frac{1}{r}=\frac{1-\theta}{2}+\frac{\theta}{1}\;\;\;\;,\;\;\;\;
\widetilde{w}=w^{-r\theta}.$$ 
So the value of $r$ has to be $\frac{2}{1+\theta}$, which implies indeed that $r\leq 2$, and we have
$$(L^2,L^1(w))_{\theta,r}=L^\frac{2}{1+\theta}(w^{-\frac{2\theta}{1+\theta}}).
$$
By using the classical formula 
$$L^p(w)^*=L^{p'}(w^{-\frac{p'}{p}})$$
(see for example (7.4.15) of \cite{Vl} with $s=0$ and the definition of $L^p(w)$ corresponding to $L^p(w^\frac{1}{p})$ in our notations) we get that
$$(L^2,L^1(w))^*_{\theta,r}=L^\frac{2}{1-\theta}(w^{\frac{2\theta}{1-\theta}}).$$

Now the relation (\ref{rrr}) gives us the weighted Strichartz
estimates 
$$\|u\|_{L^\frac{4}{n\theta}\left([0,T],\,
L^{\frac{2}{1-\theta}}(w^{\frac{2\theta}{1-\theta}})\right)}\leq 
c\|u_0\|_{L^2},$$
for all $\theta$ in $]0,\frac{2}{n}]$ satisfying $(\theta,n)\neq (1,2)$. Let us notice that we are able to include the value $\theta=0$, since in this case the estimate corresponds to the mass conservation.

By denoting 
$$p=\frac{4}{n\theta}\,\,\,,\,\,\,q=\frac{2}{1-\theta},$$
we have that 
$$\frac{2\theta}{1-\theta}=-2+\frac{2}{1-\theta}=q-2.$$
Also, the couple $(p,q)$ satisfies $\frac{2}{p}+\frac{n}{q}=\frac{n}{2}$, $(p,q,n)\neq (2,\infty,2)$ and $2\leq p,q$, therefore the weighted Strichartz estimates (\ref{Strw}) are found and the Theorem \ref{disp} is completely proved.

\subsection{Proof of the Remark \ref{dispsupp}}
$$ $$
\par

Let us choose a point $\Omega$ outside the support $A$ of the initial
data, and a point $\Omega'$ in $A$. 
Let us denote $M$ a point of the intersection of the boundary of $A$
with the part of geodesic relying 
$\Omega$ and $\Omega'$. Then
$$d(\Omega,\Omega')=d(\Omega,M)+d(M,\Omega'),$$
and 
$$d(\Omega,\Omega')\geq d(\Omega,^cA)+d(\Omega',A),$$
As the function $\frac{r}{\sinh r}$ is a decreasing function,
$$\frac{d(\Omega,\Omega')}{\sinh d(\Omega,\Omega')}\leq 
\frac{d(\Omega,^cA)+d(\Omega',A)}{\sinh \left(d(\Omega,^cA)+d(\Omega',A)\right)}.$$
By using that
$$\frac{a+b}{\sinh (a+b)}\leq c 
\frac{a}{\sinh a}\frac{b}{\sinh b},$$
we finally obtain
$$\frac{d(\Omega,\Omega')}{\sinh d(\Omega,\Omega')}\leq 
c\,\frac{d(\Omega,^cA)}{\sinh d(\Omega,^cA)}\,
\frac{d(\Omega',A)}{\sinh d(\Omega',A)}.$$
By using now the inequality (\ref{disph}) we obtain the estimate of 
the Remark \ref{dispsupp}.

\subsection{Proof of the Remark \ref{disptime}}
$$ $$
\par

We shall follow the same approach as in \S 3.2 for proving the
weighted dispersion estimate (\ref{dispw}).\\
In dimension $3$, by the explicit form (\ref{odd}), the relation
(\ref{ppp}) becomes
$$u(t,\Omega)=\frac{c}{|t|^\frac{3}{2}}
\int_{\mathbb{H}^n}u_0(\Omega')\frac{\rho}{\sinh\rho}
e^{i\frac{\rho^2}{4t}}\,d\Omega',$$
and the kernel $K$ given by
$$u(t,\Omega)=\frac{c}{|t|^\frac{3}{2}}
\int_0^\infty |u_0(r')|K(t,r,r')\sinh^{n-1}r'\,dr',$$
is exactly
$$K(t,r,r')=\frac{c}{\sinh r\sinh r'}\int_{|r-r'|}^{r+r'}e^{i\frac{y^2}{4t}} y\,dy=t\,\frac{c}{\sinh r\sinh r'}e^{i\frac{r^2+r'^2}{4t}}\cos\frac{rr'}{2t}.$$
Therefore
$$|u(t,\Omega)|\,\sinh d(\Omega,0)\leq 
\frac{c}{|t|^\frac{1}{2}}
\int_{\mathbb{H}^n}|u_0(\Omega')|\,\frac{d\Omega'}{\sinh d(\Omega',0)},$$
so the dispersion estimate of the Remark
\ref{disptime} is proved.\\


\section{Global existence and blow-up solutions}

\subsection{Global existence}
$$ $$

The Sobolev embeddings have their analogue on hyperbolic space (\cite{He})
$$\|v\|_{2^*}\leq K(n,2) \|\nabla v\|_2-\omega_n^{-\frac{2}{n}}\|v\|_2,$$
where $K(n,2)$ is the best constant for the Sobolev embeddings on $\RR^n$, 
and $\omega_n$ is the volume of the sphere $\mathbb{S}^n$. 
By interpolation between the $L^2$ and the $L^{2^*}$ norms, 
we get the Gagliardo-Nirenberg inequality for functions on $\mathbb{H}^n$
$$\|v\|_{p+1}^{p+1}\leq C_{p+1} \|v\|_2^{2+(p-1)\frac{2-n}{2}}
\|\nabla v\|_2^{(p-1)\frac{n}{2}}.$$
This inequality implies that the energy of the solution $u$ of the 
equation $(S)$ is bounded from below by
$$\|\nabla u\|_2^2\left(\frac{1}{2}-
\frac{C_{p+1}}{p+1}\|u\|_2^{2+(p-1)\frac{2-n}{2}}
\|\nabla u\|_2^{(p-1)\frac{n}{2}-2}\right).$$
As a consequence, if $p<1+\frac{4}{n}$, since the mass is conserved,
the gradient of $u$ is controlled by the energy.
Therefore the solution does not blow up and global existence occurs.\\
In the case $p=1+\frac{4}{n}$, if the mass of the initial
condition is small enough so that
$$\|u\|_2^{\frac{4}{n}}<
\frac{2+\frac{4}{n}}{2C_{2+\frac{4}{n}}},$$
then the energy controls the gradient and again, the global existence is
proved for the equation $(S)$.\\

\begin{remark}
The best constant in the Gagliardo-Niremberg inequality can be proved to be 
larger or equal to the one on $\RR^n$, but it is not obvious that it is exactly equal to it.
\end{remark}

\subsection{Blow-up solutions}
$$ $$

The power $p=1+\frac{4}{n}$ shall be proved to be the critical power, 
in the sense that the
nonlinearity is strong enough to generate solutions blowing up in a
finite time. In the following we shall show the existence of blow-up 
solutions, by analyzing an appropriate virial function on hyperbolic space.

Let $u$ be a radial solution of $(S)$ and let $h$ be a radial 
$\mathcal{C}_0^\infty(\mathbb{H}^n)$ function.
Then, by using the fact that $u$ satisfies $(S)$, we obtain the first 
derivative in time of a virial-type function
$$\partial_t\int_{\mathbb{H}^n}|u(t)|^2h \,d\Omega=
2\int_{\mathbb{H}^n} \Re \left(u(t)\overline{u}_t(t)\right)h \,d\Omega=
2\int_{\mathbb{H}^n} \Im \left(u(t)\Delta\overline{u}(t)\right)h \,d\Omega.$$
By integrating by parts
\begin{equation}\label{g'}
\partial_t\int_{\mathbb{H}^n}|u(t)|^2h \,d\Omega=
-2\int_{\mathbb{H}^n} \Im \left(u(t)\nabla\overline{u}(t)\right)\nabla h \,d\Omega.
\end{equation}
By using again the equation $(S)$ we obtain
$$\DP{t}^2\int_{\mathbb{H}^n}\mo{u}^2h=
-2\int_{\mathbb{H}^n} \Re \left((\Delta u+|u|^{p-1}u)\nabla\overline{u}\right)
\nabla h+
2\int_{\mathbb{H}^n}\Re \left(u\nabla(\Delta \overline{u}+|u|^2\overline{u})
\right)
\nabla h$$

$$=\int_{\mathbb{H}^n}
2|\nabla u|^2\Delta h
-\left(2-\frac{4}{p+1}\right)
|u|^{p+1}\Delta h-
4\Re\left(\Delta u\nabla\overline{u}\right)\nabla h-
|u|^2\Delta^2h.$$
Since $u$ and $h$ are radial functions, it follows that
$$\DP{t}^2\int_{\mathbb{H}^n}\mo{u}^2h=
\int_{\mathbb{H}^n}
4|\nabla u|^2\partial_r^2 h-
2\frac{p-1}{p+1}
|u|^{p+1}\Delta h-
|u|^2\Delta^2h.$$
By performing a density argument in the spirit of (\cite{Ca}, Lemmas 6.4.3-6), the
previous identities are valid for
$$h(\Omega)=h(\cosh r,\sinh r\omega)=r^2.$$
By direct computation we obtain 
$$\partial_r^2 r^2=2,$$
and
$$\Delta r^2=2+2(n-1)\frac{\cosh r}{\sinh r}r.$$
Therefore we can display the energy in the second derivative of the 
virial 
\begin{equation}\label{vir}
\DP{t}^2\int_{\mathbb{H}^n}\mo{u}^2r^2=
16E(u)-\int_{\mathbb{H}^n}|u|^2\Delta^2r^2-
2\frac{p-1}{p+1}\int_{\mathbb{H}^n} |u|^{p+1}2(n-1)
\left(\frac{\cosh r}{\sinh r}r-1\right)-
\end{equation}

$$-\left(-\frac{16}{p+1}+
2\frac{p-1}{p+1}2n\right)\int_{\mathbb{H}^n} |u|^{p+1}.$$
We have
\begin{equation}\label{delta2}
\Delta^2 r^2=4(n-1)\frac{1}{\sinh^3r}(r\cosh r-\sinh r)+
2(n-1)^2\frac{\cosh r}{\sinh^3r}(\cosh r\sinh r-r),\end{equation}
nonnegative for all $r\geq 0$, since 
$$\left\{\begin{array}{c}r\cosh r-\sinh r\geq 0\\
\cosh r\sinh r-r\geq 0\end{array}\right..$$
Therefore, since $p$ is considered larger than $1$, the second and the third 
term in the right hand side of (\ref{vir}) are negative. The last one is also negative, provided that 
$$-\frac{16}{p+1}+
4n\frac{p-1}{p+1}\geq 0,$$
that is exactly
$$p\geq 1+\frac{4}{n}.$$
Let us notice, 
by (\ref{delta2}), that $\Delta^2 r^2$ is bounded between two 
positive constants depending only on the dimension $n$,
$$0<k_n<\Delta^2 r^2<K_n.$$

In conclusion, if $u$ is a solution of $(S)$ of initial data satisfying
$$16E(u_0)<k_n\|u_0\|^2_2,$$
there exists a constant $C$ such that for any $t$
\begin{equation}\label{negder}
\DP{t}^2\int_{\mathbb{H}^n}\mo{u(t)}^2r^2<C<0.\end{equation}
It follows that there is a finite time $T$ for which 
$$\lim_{t\tend T}\int_{\mathbb{H}^n}\mo{u(t)}^2r^2=0.$$
Then, using the uncertainty principle
$$\left(\int_{\mathbb{H}^n}\mo{u}^2\right)^2\leq c
\left(\int_{\mathbb{H}^n}\mo{u}^2r^2\right)
\left(\int_{\mathbb{H}^n}\mo{\nabla u}^2\right),$$
tells us that
$$\lim_{t\tend T}\int_{\mathbb{H}^n}\mo{\nabla u(t)}^2=+\infty,$$
so the Theorem \ref{blup} is proved.

Finally, let us mention that in the Euclidean case, one has
$$\DP{t}^2\int_{\mathbb{R}^n}\mo{u}^2|x|^2=16E(u),$$
so the argument used before does not work for any solutions of null energy. \\

\section{Appendix A}
In this Appendix the following Proposition shall be proved.\\

{\bf{Proposition \ref{F}}}
{\it{For all $m\geq 1$ integer, the following identity holds 
$$(\ref{derexp})\,\,\,
\left(\frac{\partial_s}{\sinh s}\right)^me^{i\frac{s^2}{4t}}=
\sum_{k=1}^m \frac{e^{i\frac{s^2}{4t}}}{t^k}\,F_k^m(s),$$
where, modulo constants,  
$$(\ref{idFm})\,\,\,
F_m^m(s)=\left(\frac{s}{\sinh s}\right)^m,$$
and, for $1\leq k<m$, $F_k^m(s)$ can be written
$$(\ref{idF})\,\,\,
F_k^m(s)=\sum_{r=0}^k\,\,\sum_{\left\{\begin{array}{c}i_1+\cdot\cdot\cdot+i_r=m-k\\i_j\geq 1\end{array}\right.}\left(\frac{s}{\sinh s}\right)^{k-r}
\left(\frac{\partial_s}{\sinh s}\right)^{i_1}
\frac{s}{\sinh s}...
\left(\frac{\partial_s}{\sinh s}\right)^{i_r}
\frac{s}{\sinh s}.$$
Moreover, for $\alpha\in\{0,1\}$ and $1\leq k\leq m$
$$(\ref{estFk})\,\,\,
\left|\partial_s^\alpha\left(\frac{\sinh s}{s}\,F_k^m(s)\right)\right|\leq 
c\,s^\alpha\left(\frac{s}{\sinh s}\right)^{m-1},$$
and the second derivative of $\frac{\sinh s}{s}\,F_k^m(s)$ is bounded for $m\geq 2$.}}\\

\begin{proof}
Each time we differentiate the exponential $e^{i\frac{s^2}{4t}}$ we obtain a term with one negative power of $t$ times $\frac{s}{\sinh s}$. 
Therefore the terms $F_k^m(s)$ are the ones corresponding to the derivation 
of $e^{i\frac{s^2}{4t}}$ $k$ times, and they contain $k$ times 
the function $\frac{s}{\sinh s}$ or its derivatives.\\

\noindent
We have denoted in (\ref{idF}) by $k-r$ the number of times 
$\frac{s}{\sinh s}$ appears, where $0\leq r\leq k$ ~; the other $r$ terms are 
derivatives of $\frac{s}{\sinh s}$. 
Since only $k$ derivatives have fallen on the exponential term, there are 
$m-k$ left derivatives which act on the remainder term. 
That is the reason for which in (\ref{idF}) we have the sum
$$i_1+...+i_r=m-k.$$

For showing the estimates of $F_k^m(s)$ we shall need the following Lemma.

\begin{lemma}\label{Pk}
The following property  
\begin{equation}\label{derex}\exists\; c_l>0 \;\;\mbox{ such that }\;\;
\left|\left(\frac{\partial_\rho}{\sinh\rho}\right)^le^{i\rho^2}\right|\leq 
c_l\left(\frac{\rho}{\sinh\rho}\right)^l, \forall \rho\geq 0,\end{equation}
is valid for all positive integer $l$. From it follows also that
\begin{equation}\label{dersh} P(l) \,\,:\,\, \,\,\exists\;c_l>0 \;\;\mbox{ such that }\;\;
\left|\left(\frac{\partial_\rho}{\sinh\rho}\right)^l
\frac{\rho}{\sinh\rho}\right|\leq 
c_l\left(\frac{\rho}{\sinh\rho}\right)^{l+1}, \forall \rho\geq 0.\end{equation}
\end{lemma}
\begin{proof}
We shall split the proof into two cases, one when $\rho$ is small, and one when it is large.\\

{\underline{The case $\rho<1$}}.\par

The function
$$G_l(\rho)=\left(\frac{\partial_\rho}{\sinh\rho}\right)^le^{i\rho^2}$$
is a $\mathcal{C}^\infty$ even function. 
It is obvious for $l=0$, and for other values it can be proved by 
induction as follows. 
Supposing that the fact is true for $l$, we have the existence of a $\mathcal{C}^\infty$ 
odd function $H_l$ such that
$$G'_l(\rho)=H_l(\rho).$$
Therefore
$$G_{l+1}(\rho)=\frac{\partial_\rho}{\sinh\rho}G_l(\rho)=\frac{G'_l(\rho)}{\sinh \rho}=\frac{H_l(\rho)}{\sinh\rho},$$
and the induction argument is done.\par
As $\rho$ is small, and $G_l$ is a $\mathcal{C}^\infty$ 
even fuction, we obtain the existence of a constant $c$ such that
$$|G_l(\rho)|\leq c.$$
Again because $\rho$ is small, we have the quotient $\frac{\rho}{\sinh\rho}$ bounded, and 
we can find for all integer $l$ a constant $c_l$ such that
$$|G_l(\rho)|\leq c_l\left(\frac{\rho}{\sinh\rho}\right)^l.$$
So in the case $\rho<1$ the property $(\ref{derex})$ is satisfied.\\

{\underline{The case $\rho\geq 1$}}.\par

Let us perform the change of variable $r=\cosh\rho$. Then 
$$\frac{\partial_\rho}{\sinh\rho}=\partial_r,$$
and $\rho$ can be written
$$\rho=\log (r+\sqrt{r^2-1}).$$
With these notations, the property $(\ref{derex})$ becomes
$$\left|\partial_r^le^{i\left(\log (r+\sqrt{r^2-1})\right)^2}\right|\leq 
c_l\left(\frac{\log (r+\sqrt{r^2-1})}{\sqrt{r^2-1}}\right)^l, \forall \rho\geq 1.$$

Let us denote
$$a(r)=r+\sqrt{r^2-1},$$
and $$b(r)=\log a(r).$$

Since $r$ is larger than $\cosh 1$, then for proving the property $(\ref{derex})$ it suffices to prove 
that
\begin{equation}\label{expib}
\left|\partial_r^le^{ib^2(r)}\right|\leq 
c_l\frac{b^l(r)}{r^l}, \forall \,r\geq 3.\end{equation}

In the following we shall use the following general formula
\begin{equation}\label{Faa}
\partial_r^lf(g(r))=\sum_{q=1}^lc_qf^{(q)}(g(r))\sum_{p_1+...+p_q=l}c_{p_1,..,p_q}\partial_r^{p_1}g(r)...\partial_r^{p_q}g(r).\end{equation}
This is a weaker form of the Faa di Bruno formula, in which the constants are described. 

Let us notice first that
 \begin{equation}\label{a}
|\partial_r^pa(r)|\leq cr^{1-p},
\end{equation}
and that, since $r$ is large enough,
\begin{equation}\label{loga}
|\log^{(q)}(a(r))|\leq cr^{-q},\end{equation}
for all $q\geq 1$.
Then, by using (\ref{Faa}), for $l\geq 1$,
$$\partial_r^l\log(a(r))=\sum_{q=1}^lc_q\log^{(q)}(a(r))\sum_{p_1+...+p_q=l}c_{p_1,..,p_q}\partial_r^{p_1}a(r)...\partial_r^{p_q}a(r).$$
Now the estimates (\ref{a}) and (\ref{loga}) tell us that
$$|\partial_r^l\log(a(r))|\leq c\sum_{q=1}^lr^{-q}\sum_{p_1+...+p_q=l}r^{1-p_1}...r^{1-p_q}=cr^{-l}.$$
Therefore we have obtained that for $l\geq 1$,
\begin{equation}\label{b}
|\partial_r^l b(r)|\leq cr^{-l}.\end{equation}
Let us compute now the derivatives of $b^2$
$$\partial_r^l b^2(r)=\sum_{q=0}^l\partial_r^q b(r)\partial_r^{l-q} b(r)=2b(r)\partial_r^l b(r)+\sum_{q=1}^{l-1}\partial_r^q b(r)\partial_r^{l-q} b(r).$$
By using (\ref{b}) we get
$$|\partial_r^l b^2(r)|\leq c|b(r)|r^{-l}+\widetilde{c}\sum_{q=1}^{l-1}r^{-q}r^{-(l-q)}\leq (c|b(r)|+c)r^{-l}.$$
As the variable $r$ is large, it follows that
$$b(r)=\log (r+\sqrt{r^2-1})\geq c,$$
and therefore
\begin{equation}\label{b2}
|\partial_r^l b^2(r)|\leq c|b(r)|r^{-l}.\end{equation}

Finally, by using again the general formula (\ref{Faa}), 
$$\left|\partial_r^le^{ib^2(r)}\right|\leq c\sum_{q=1}^l\sum_{p_1+...+p_q=l}\partial_r^{p_1}b^2(r)...\partial_r^{p_q}b^2(r),$$
so the estimates (\ref{b2}) allows us to write
$$\left|\partial_r^le^{ib^2(r)}\right|\leq c\sum_{q=1}^lb^q(r)r^{-l}.$$
We have already noticed that $b(r)$ is large, therefore
$$\left|\partial_r^le^{ib^2(r)}\right|\leq cb^l(r)r^{-l},$$
and we finally proved (\ref{expib}), which implies the wanted property $(\ref{derex})$.\\

Property (\ref{dersh}) can be proved by induction. $P(0)$ is obvious. Suppose that the properties $P(1)$,...,$P(l-1)$ are true. Then the development (\ref{derexp}) for $m=l+1$ will contain, by (\ref{idFm}) and (\ref{idF}) only one term of the type 
$$\left(\frac{\partial_\rho}{\sinh\rho}\right)^{l}
\frac{\rho}{\sinh\rho},$$
all others being of order less than $l$. More precisely, we have, modulo constants,
$$\left(\frac{\partial_s}{\sinh s}\right)^{l+1}e^{i\frac{s^2}{4t}}=
\sum_{k=1}^{l+1} \frac{e^{i\frac{s^2}{4t}}}{t^k}\,F_k^{l+1}(s)=$$
$$=\frac{e^{i\frac{s^2}{4t}}}{t}\left(\frac{\partial_\rho}{\sinh\rho}\right)^{l}
\frac{\rho}{\sinh\rho}+\frac{e^{i\frac{s^2}{4t}}}{t^{l+1}}\left(\frac{s}{\sinh s}\right)^{l+1}+$$
$$+\sum_{k=2}^l\frac{e^{i\frac{s^2}{4t}}}{t^k}\sum_{r=0}^k\,\,\sum_{\left\{\begin{array}{c}i_1+\cdot\cdot\cdot+i_r=l+1-k\\i_j\geq 1\end{array}\right.}\left(\frac{s}{\sinh s}\right)^{k-r}
\left(\frac{\partial_s}{\sinh s}\right)^{i_1}
\frac{s}{\sinh s}...
\left(\frac{\partial_s}{\sinh s}\right)^{i_r}
\frac{s}{\sinh s}.$$
In the last term all index $i_j$ are less than $l$, and the properties $P(1)$,...,$P(l-1)$ tell us that the term is upper-bounded by 
$$\left(\frac{s}{\sinh s}\right)^{l+1}.$$
Therefore, by using (\ref{derex}) we obtain that also the first term in the right-hand-side is bounded by the same quantity
$$\left|\left(\frac{\partial_\rho}{\sinh\rho}\right)^{l+1}
\frac{\rho}{\sinh\rho}\right|\leq 
c_k\left(\frac{\rho}{\sinh\rho}\right)^{l+1},$$
so the property $P(l)$ is proven and the induction argument is complete.

\end{proof}

Let us return now to the estimates of Proposition \ref{F}. \\

For $k=m$, estimate (\ref{estFk}) when $\alpha =0$ is obvious using (\ref{idFm}). For $0\leq k<m$, by (\ref{idF}) and we get that
$$\left|\frac{\sinh s}{s}
F_k^m(s)\right|\leq\sum_{r=0}^k\,\,\sum_{\left\{\begin{array}{c}i_1+\cdot\cdot\cdot+i_r=m-k\\i_j\geq 1\end{array}\right.}\left(\frac{s}{\sinh s}\right)^{k-r-1}
\left(\frac{s}{\sinh s}\right)^{i_1+1}...
\left(\frac{s}{\sinh s}\right)^{i_r+1},$$
so estimate (\ref{estFk}) when $\alpha =0$ follows by using (\ref{dersh}).\\

We start now to look at the derivatives. One has
\begin{equation}\label{1ereder}\partial_s\left(\frac{\sinh s}{s}\,F_k^m(s)\right)=\frac{s\cosh s-\sinh s}{s^2}F_k^m(s)+\frac{\sinh^2 s}{s}\frac{\partial_s}{\sinh s}F_k^m(s)\end{equation}
The forms (\ref{idFm}) and (\ref{idF}) of $F_k^m(s)$, together with estimate (\ref{dersh}), imply
$$\left|\partial_s\left(\frac{\sinh s}{s}\,F_k^m(s)\right)\right|\leq \frac{s\cosh s-\sinh s}{s^2}\left(\frac{s}{\sinh s}\right)^m+\frac{\sinh^2 s}{s}\left(\frac{s}{\sinh s}\right)^{m+1}.$$
We deduce that
$$\left|\partial_s\left(\frac{\sinh s}{s}\,F_k^m(s)\right)\right|\leq s\left(\frac{s}{\sinh s}\right)^{m-1}\left(\frac{s\cosh s-\sinh s}{s^2\sinh s}+1\right),$$
and estimate (\ref{estFk}) for $\alpha =1$ follows.\\

Knowing (\ref{1ereder}) one can get the expression of the second derivative
$$\partial_s^2\left(\frac{\sinh s}{s}\,F_k^m(s)\right)=\frac{s^2\sinh s-2(s\cosh s-\sinh s)}{s^3}F_k^m(s)+$$
$$+\frac{3s\cosh s\sinh s-2\sinh^2 s}{s^2}\frac{\partial_s}{\sinh s}F_k^m(s)+
\frac{\sinh^3 s}{s}\left(\frac{\partial_s}{\sinh s}\right)^2F_k^m(s).$$
By using again the forms (\ref{idFm}) and (\ref{idF}) of $F_k^m(s)$ and the estimate (\ref{dersh}) 
$$\left|\partial_s^2\left(\frac{\sinh s}{s}\,F_k^m(s)\right)\right|\leq \left(\frac{s}{\sinh s}\right)^{m-2}\left(\frac{s^2\sinh s-2(s\cosh s-\sinh s)}{s\sinh^2 s}+\right.$$
$$\left.+s\frac{3s\cosh s-2\sinh s}{\sinh^2 s}+\frac{s^3}{\sinh s}\right).$$
As the term in the brackets of the right-hand-side is bounded, for $m\geq 2$ we obtain that the second derivative of $\frac{\sinh s}{s}\,F_k^m(s)$ is bounded, and Proposition \ref{F} is completely proved.
\end{proof}


\begin{thebibliography}{99999}
\bibitem{AG} S. Alinhac, P. G\'erard,
Op\'erateurs pseudo-diff\'erentiels et th\'eor\`eme de Nash-Moser, 
Paris InterEditions CNRS 1991.

\bibitem{An} J.-P. Anker, L. Ji, 
Heat kernel and Green function estimates on noncompact symmetric
spaces, 
Geom. and Funct. Anal. 9, (1999), no. 6, 1035-1091.

\bibitem{Val2} V. Banica,
On the nonlinear Schr\"odinger dynamics on $\mathbb{S}^2$,
J. Math. Pures Appl. 83 (2003), no. 1, 77-98.

\bibitem{B} H. Bateman, 
Higher transcendental functions, 
New York NY Toronto London McGraw-Hill 1953-1955.

\bibitem{BL} J. Bergh, J. L\"ofstr\"om,
Interpolation spaces an introduction, 
Berlin Heidelberg New York NY Springer 1976.

\bibitem{Bo} J. Bourgain,
Fourier transformation restriction phenomena for certain lattice subsets
and application to the nonlinear evolution equations I. Schr\"odinger
equations,
Geom. and Funct. Anal. 3 (1993), no. 2, 107-156.

\bibitem{BGT12} N. Burq, P. G\'erard, N. Tzvetkov,
An instability property of the nonlinear Schr\"odinger equation on
$\mathbb{S}^d$,
Math. Res. Lett. 9 (2002), no. 2-3, 323-335.

\bibitem{BGT0} N. Burq, P. G\'erard, N. Tzvetkov,
Strichartz inequalities and the nonlinear Schr\"odinger equation on compact manifolds, 
Amer. J. Math. 126 (2004), no. 3, 569-605.

\bibitem{BGT32} N. Burq, P. G{\'e}rard, N. Tzvetkov,
In{\'e}galit{\'e}s de Sogge bilin{\'e}aires et {\'e}quation de
Schr{\"o}dinger non-lin{\'e}aire,
S{\'e}minaire {\'E}quations aux D{\'e}riv{\'e}es partielles,
{\'E}cole Polytechnique, Palaiseau, mars 2003.


\bibitem{Ca} T. Cazenave, 
An introduction to nonlinear Schr{\"o}dinger equations, 
Textos de M{\'e}todos Matem{\'a}ticos 26, Instituto de 
Matem{\'a}tica-UFRJ, Rio de Janeiro, RJ 1996. 


\bibitem{CW} T. Cazenave, F. Weissler,
The Cauchy problem for the nonlinear Schr{\"o}dinger equation
in $\HH^1$,
Manuscripta Math. 61 (1988), no. 4, 477-494.


\bibitem{CCT} M. Christ, J. Colliander, T. Tao,
Asymptotics, frequency modulation, and low regularity ill-posedness for
canonical defocusing equations,
Amer. J. Math. 125 (2003), no. 6, 1235-1293.

\bibitem{D} E. B. Davies, 
Heat kernels and spectral theory, 
Cambridge Univ. Press 1989.

\bibitem{DM} E. B. Davies, N. Mandouvalos, 
Heat kernel bounds on hyperbolic space and Kleinian groups, 
Proc. London Math. Soc. 57 (1988), no. 3, 182-208.

\bibitem{Vl} V. Georgiev,
Semilinear hyperbolic equations,
Mathematical Society of Japan, Tokyo 2000.

\bibitem{GV} J. Ginibre, G. Velo,
The global Cauchy problem for the nonlinear Schr{\"o}dinger equation,
Ann. I. H. P. Analyse non-lin{\'e}aire 2 (1985), no. 4, 309-327.

\bibitem{Gl} R. T. Glassey,
On the blowing up of solutions to the Cauchy problem for nonlinear
Schr{\"o}dinger equations,
J. Math. Phys. 18 (1977), 1794-1797.

\bibitem{He} E. Hebey, 
Nonlinear analysis on manifolds: Sobolev spaces and inequalities, 
NY Courant Institute of Mathematical Sciences, New York 1999.

\bibitem{Hel} S. Helgason,
Geometric analysis on symmetric spaces,
Publication Providence RI American mathematical society (1994).

\bibitem{Ka} O. Kavian, 
A remark on the blowing-up of solutions to the Cauchy problem 
for nonlinear Schr{\"o}dinger equations, 
Trans. Amer. Math. Soc. 299 (1987), no. 1, 193-203.

\bibitem{KT} M. Keel, T. Tao, 
Endpoint Strichartz estimates, 
Amer. J. Math. 120 (1998), no. 5, 955-980.

\bibitem{Me} F. Merle,
Determination of blow-up solutions with minimal mass for nonlinear
Schr{\"o}dinger equation with critical power,
Duke Math. J. 69 (1993), 427-454.

\bibitem{Vi} V. Pierfelice, Weighted Strichartz estimates for the radial perturbed 
Schr\"odinger equation on the hyperbolic space, 
Manuscripta Mathematica 120 (2006), 377-389.

\bibitem{Ta} D. Tataru,
Strichartz estimates in the hyperbolic space and global existence for
the semilinear wave equation,
Trans. Amer. Math. Soc. 353 (2001), no. 2, 795-807.

\bibitem{Ter} A. A. Terras, 
Harmonic analysis on symmetric spaces and applications 1, 
New York NY Berlin Heidelberg Springer 1985.

\bibitem{To} P. Tomas,
A restriction theorem for the Fourier transform,
Bull. Amer. Math. Soc. 81 (1975), 177-178.

\bibitem{Ya} K. Yajima,
Existence of solutions for Schr{\"o}dinger evolution equations,
Comm. Math. Phys. 110 (1987), no. 3, 415-426.

\bibitem{Za} V. E. Zakharov,
Collapse of Lagmuir waves,
Sov. Phys. JETP 35 (1972), 908-914.


\end{thebibliography}
\end{document}